\documentclass[10pt]{article}
\usepackage{amsmath, amsthm,  latexsym, amssymb, amsfonts, epsfig,
caption2, graphicx}
\usepackage[dvipdfm]{hyperref}
\usepackage{mathrsfs}

\def\cajita{\rule{5pt}{5pt}}

\renewcommand{\d}{\mathrm{d}}
\newcommand{\esssup}{\mathop{\mathrm{esssup}}}
\newcommand{\ep}{\varepsilon}

\newcommand{\dw}{\downarrow}

\newcommand{\wt}{\widetilde}
\newcommand{\pr}{\prime}
\newcommand{\N}{\mathbb{N}}
\newcommand{\R}{\mathbb{R}}

\newcommand{\E}{\mathbb{E}}
\newcommand{\PP}{\mathbb{P}}

\newcommand{\HH}{\mathcal{H}}
\newcommand{\F}{\mathcal{F}}
\newcommand{\lt}{\left}
\newcommand{\rt}{\right}
\newtheorem{theorem}{Theorem}[section]
\newtheorem{proposition}[theorem]
{Proposition}
\newtheorem{lemma}[theorem]{Lemma}
\par
\newtheorem{definition}[theorem]
{Definition}
\par

\par
\par
\newtheorem{remark}[theorem]{Remark}
\par

\par
\newtheorem{example}[theorem]
{Example}
\par

\begin{document}

\title{Nonlinear Fractional Backward
Doubly Stochastic Differential Equations with Hurst Parameter in
(1/2,1)}
\author{ Shuai Jing\thanks{This work is
supported by European Marie Curie Initial Training Network (ITN)
project: $"$Deterministic and Stochastic Controlled Systems and
Application$"$, FP7-PEOPLE-2007-1-1-ITN, No. 213841-2, the National
Basic Research  Program of China (973 Program) grant No.
2007CB814900 (Financial Risk) and the NSF of China (No.11071144).}\\
\small{ D\'epartement de Math\'ematiques, Universit\'e de Bretagne
Occidentale,} \small{29285 Brest C\'edex, France} \\
\small{School of Mathematics, Shandong University, 250100, Jinan,
China}\\\small{E-mail:
\href{mailto:shuai.jing@univ-brest.fr}{shuai.jing@univ-brest.fr}}}
\maketitle

\begin{abstract}{We first state a special type of It\^o formula
involving stochastic integrals of both standard and fractional
Brownian motions. Then we use Doss-Sussman transformation to
establish the link between backward doubly stochastic differential
equations, driven by both standard and fractional Brownian motions,
and backward stochastic differential equations, driven only by
standard Brownian motions. Following the same technique, we further
study associated nonlinear stochastic partial differential equations
driven by fractional Brownian motions and partial differential
equations with stochastic coefficients.}
\end{abstract}

\bigskip


\section{Introduction}\label{sec:1}
The theory of backward stochastic differential equations (BSDEs) and
that of fractional Brownian motion (fBm) had  developed
simultaneously in their own separated directions for many years
until Bender \cite{Be} gave an explicit solution to a linear BSDE
driven by fBm in 2005. In 2009 Hu and Peng  \cite{HP} stated a more
general theory on fractional BSDEs by using the so-called
quasi-expectation, but their case is still limited. The
non-semimartingale property of the fBm (except the case of Hurst
parameter $H=1/2$, where it becomes a Brownian motion) makes it
thorny to handle. Being not a semimartingale means there is no
martingale representation theory for the fBm, which is crucial in
the general BSDE theory (see the pioneering work of Pardoux and Peng
\cite{PP3}). In Jing and Le\'on \cite{Jl}, we tried to combine the
fBm and BSDEs in another way: We transformed a semilinear backward
doubly stochastic differential equation (BDSDE) driven by both a
standard and a fractional Brownian motions with $H\in (0,1/2)$ into
a BSDE without integral with respect to the fBm, which turns out
easier to deal with. The integral w.r.t. the fBm in the BDSDE was
interpreted in the sense of the extended divergence operator and the
integral was supposed to be linear w.r.t. the solution process. This
allowed to apply the efficient tool of nonanticipating Girsanov
transformation, developed in Buckdahn \cite{bu}. However, this
method is restricted to semilinear BDSDEs.

In this paper, we deal with BDSDEs for which the integrand of the
integral with respect to the fBm is not necessarily linear with the
solution process, and the Hurst coefficient $H$ is supposed to
belong t o the interval$(1/2,1)$. Unlike the more irregular case
$H<1/2$, the stochastic integrals with respect to an fBm with
$H>1/2$ can be defined in different ways. So they can be defined
with the help of the divergence operator in the frame of the
Malliavin calculus, see Decreusefont and \"Ust\"unel \cite{Du} and
Al\`{o}s {\it et al.} \cite{Amn} (Notice that the Wick-It\^o
integral defined in Duncan {\it et al.} \cite{Dhp} coincides with
the first one). They can also be defined pathwise as generalized
Riemann-Stieltjes integral (see Z\"ahle \cite{Za1} and \cite{Za2})
or with the help of the rough path theory (see Coutin and Qian
\cite{Cq}). For a complete list of references we refer to the two
books by Biagini {\it et al.} \cite{BHOZ} and Mishura \cite{Mi}.

Our approach to BDSDEs with an fBM is inspired by the work of
Buckdahn and Ma \cite{BM}. In their study of stochastic PDEs driven
by a Brownian motion $B$ the authors of \cite{BM} used  BDSDEs
driven by $B$ as well as an independent Brownian motion; the
integral with respect to $B$ is interpreted in Stratonovich sense.
This allowed the application of the Doss-Sussman transformation in
order to transform the BDSDE into a BSDE without integral with
respect to $B$. On the other hand, the pathwise integral with
respect to the fBm plays a role which is comparable with that of the
Stratonovich integral in the classical theory. Nualart and
R\u{a}\c{s}canu \cite{Nr} used the pathwise integral to solve
(forward) stochastic differential equations driven by an fBm. For
some technical reasons (such as the lack of H\"older continuity, see
Remark \ref{rmk2}), we shall make use of the Russo-Vallois integral
developed by Russo and Vallois in a series of papers (\cite{Rv1},
\cite{Rv1}, \cite{Rv1}, {\it etc}.). Under standard assumptions
which allow to apply the Doss-Sussman transformation, we associate
the BDSDE driven by both a standard Brownian motion $W$ and an fBm
$B$,
\begin{equation}\label{eq:1.2}
U_s^{t,x}=\Phi(X_0^{t,x})+\int^s_0f(r,X_r^{t,x},
U_r^{t,x},V_r^{t,x})\d r + \int^s_0 g(U_r^{t,x})\d
B_r-\int^s_0V_r^{t,x}\dw\d W_r, \ s\in[0,t], \footnote{$\int^t_0
\cdot\dw \d W_s$ \textrm{indicates that the integral is considered
as the It\^o backward one.}}
\end{equation}
with the BSDE driven only by the Brownian motion $W$,
\begin{equation} \label{eq:1.1}
Y_s^{t,x}=\Phi(X_0^{t,x})+\int^s_0\tilde{f}(r,X_r^{t,x},Y_r^{t,x},Z_r^{t,x})\d
r -\int^s_0 Z_r^{t,x}\dw \d W_r, \ s\in[0,t].
\end{equation}
Here $\tilde{f}$ will be specified in Section 4;  it is a driver
with quadratic growth in $z$. We point out that the classical BDSDEs
were first studied by Pardoux and Peng \cite{PP2} and our BSDE
(\ref{eq:1.1}) is a quadratic growth BSDE, which was studied first
by Kobylanski \cite{Ko}. In the works of Pardoux and Peng \cite{PP2}
and Buckdahn and Ma \cite{BM}, i.e., when the Hurst parameter
$H=1/2$, one can solve the BDSDE directly and get the square
integrability of the solution process. However, in the fractional
case $(H\neq1/2)$, to our best knowledge, there does not exist a
direct way to solve the BDSDE (\ref{eq:1.2}), and as it turns out in
Theorem \ref{thm:sol}, we can only get that the conditional
expectation of $\int^t_0\lt|Z^{t,x}_r\rt|^2\d r$ is bounded by an
a.s. finite process. This is also the reason that instead of using
the space of square integrable processes, we use the space of a.s.
conditionally square integrable processes (see the definition of the
space $\HH_T^2(\R^d)$ in Section 2).

A celebrated contribution of the BSDE theory consists in giving a
form of probabilistic interpretation, nonlinear Feynman-Kac formula,
to the solutions of PDEs (see, for instance, Peng \cite{Pe}, Pardoux
and Peng \cite{PP1}). As indicated in Kobylanski \cite{Ko}, the
quadratic growth BSDE (\ref{eq:1.1}) is connected to the semilinear
parabolic PDE
\begin{equation}\label{eq:1.3}
\left\{
\begin{array}{cllc}
\frac{\partial u}{\partial t}(t,x)&=&\mathscr{L} u (t,x) -
\tilde{f}(t,x,u(t,x),\sigma(x)^T\frac{\partial}{\partial x}u(t,x)),&
(t,x)\in (0,T)\times
\R^n;\\
u(0,x)&=&\Phi(x),&x\in\R^n,
\end{array}\right.
\end{equation}
where $\mathscr{L}$ is the infinitesimal operator of a Markov
process. Hence, it is natural for us to consider the form of equation
(\ref{eq:1.3}) after Doss-Sussman transformation and we prove that
it becomes the following semilinear SPDE
\begin{equation}\label{eq:1.4}
\left\{
\begin{array}{lc}
\d u(t,x)=\left[\mathscr{L} u (t,x) -
{f}(t,x,u(t,x),\sigma(x)^T\frac{\partial}{\partial x}u(t,x))\rt]\d t
+ g(u(t,x))\d B_t,& (t,x)\in (0,T)\times
\R^n;\\
u(0,x)=\Phi(x),&x\in\R^n.
\end{array}\right.
\end{equation}

We emphasize that this paper can not be considered as a
generalization of \cite{Jl}. The reason is that, firstly, the Hurst
parameters are distinct; secondly, the stochastic integrals with
respect to the fBm are of different types.

We organize the paper as follows. In section 2, we recall some basic
facts about the fractional Brownian motion, we give the general
framework of our work, and we recall the definition of the backward
Russo-Vallois integral as well as some of its properties. In section
3 we prove a type of It\^o formula involving integrals with respect
to both standard and fractional Brownian motions, which will play an
important role in the following sections. We perform a Doss-Sussman
transformation in Section 4 to transform a nonlinear BDSDE
(\ref{eq:1.2}) into a BSDE (\ref{eq:1.1}) and show the relationship
between their solutions. In particular, we show that BSDE
(\ref{eq:1.1}) has a unique solution $(Y^{t,x},Z^{t,x})$, and the
couple of processes $(U^{t,x},V^{t,x})$ associated with
$(Y^{t,x},Z^{t,x})$ by the inverse Doss-Sussman transformation is
the unique solution of BDSDE ({\ref{eq:1.2}}). Finally, the
stochastic PDE associated with BDSDE (\ref{eq:1.2}) is briefly
discussed in Section 5.
\section{Preliminaries}
\subsection{Fractional Brownian Motion and General Setting}

In this subsection we recall some basic results on the fBm and the
related setting. For a more complete overview of the theory of fBm,
we refer the reader to Biagini {\it et al.} \cite{BHOZ} and Mishura
\cite{Mi}.

Let $(\Omega^\pr,\F^\pr, \PP^\pr)$ be a classical Wiener space with
time horizon $T>0$, i.e., $\Omega^\pr=C_0([0,T];\R)$ denotes the set
of real-valued continuous functions starting from zero at time zero,
endowed with the topology of the uniform convergence,
$\mathcal{B}(\Omega^\pr)$ is the Borel $\sigma$-algebra on
$\Omega^\pr$ and $\PP^\pr$ is the unique probability measure on
$(\Omega, \mathcal{B}(\Omega^\pr))$ with respect to which the
coordinate process $W^0_t(\omega^\pr)=\omega^\pr(t),$  $t\in[0,T],$
$\omega^\pr\in\Omega^\pr$ is a standard Brownian motion. By $\F^\pr$
we denote the completion of $\mathcal{B}(\Omega^\pr)$ by all
$\PP^\pr$-null sets in $\Omega^\pr$. Given $H\in(1/2,1)$, we define
$$
B_t=\int_0^t K_H(t,s)\d W^0_s,\quad t\in[0,T],
$$
where  $K_H$ is the kernel of the fBm with parameter $H\in(1/2,1)$:
$$
K_H(t,s)=C_Hs^{1/2 - H}\int_s^t u^{H-1/2}(u-s)^{H-3/2}\d u,
$$
with $C_H =\sqrt{\frac{H}{(2H-1)\beta (2-2H, H - {1}/{2})}}$. It is
well known that such defined process $B$ is a one-dimensional fBm,
i.e., it is a Gaussian process with zero mean and covariance
function
$$
R_H(t,s):= \mathbb{E}\left[B_tB_s\right] = \frac{1}{2}\left(t^{2H}+
s^{2H}-|t-s|^{2H}\right), \quad s, t\in [0,T].
$$

We let $\{W_t: {0\leq t\leq T}\}$ be the coordinate process on the
classical Wiener space $(\Omega^{\pr\pr},{\cal F} ^{\pr\pr},
\mathbb{P}^{\pr\pr})$ with $\Omega^{\pr\pr}=C_0([0,T]; \R^d)$, which
is a $d$-dimensional Brownian motion with respect to the Wiener
measure $\PP^{\pr\pr}$. We put $(\Omega, {\cal F}^0,
\mathbb{P})=(\Omega^{\pr},{\cal F}^{\pr}, \mathbb{P}^{\pr})
\otimes(\Omega^{\pr\pr},{\cal F} ^{\pr\pr}, \mathbb{P}^{\pr\pr})$
and let $\mathcal {F}=\mathcal{F}^0\vee \mathcal{N}$, where
$\mathcal N$ is the class of the $ \mathbb{P}$-null sets. We denote
again by $B$ and $W$ the canonical extensions of the fBm $B$ and of
the Brownian motion $W$ from $\Omega^{\pr}$ and $\Omega^{\pr\pr}$,
respectively, to $\Omega$.

We let $\mathcal{F}_{[t,T]}^W = \sigma\{W_T-W_s, t\leq s\leq
T\}\vee\cal N$, $\mathcal{F}_t^{B}=\sigma\{B_s, 0\leq s\leq
t\}\vee\cal N,$ and $\mathcal{G}_t = \mathcal{F}_{[t,T]}^W
\vee\mathcal{F}_t^{B}, t\in[0,T].$ Let us point out that
$\mathcal{F}_{[t,T]}^W$ is decreasing and $\mathcal{F}_t^{B}$ is
increasing in $t$, but $\mathcal{G}_t$ is neither decreasing nor
increasing. We denote the family of $\sigma$-fields
$\{\mathcal{G}_t\}_{0\leq t\leq T}$ by $\mathbb{G}$. Moreover, we
also introduce the backward filtrations $\mathbb{H}=\{\mathcal {H}_t
=\mathcal{F}_{[t,T]}^W \vee \mathcal{F}_T^{B}\}_ {t\in[0,T]}$ and
$\mathbb{F}^W=\{\mathcal{F}_{[t,T]}^W\}_{t\in[0,T]}$.

Finally, we denote by $C(\mathbb{H},[0,T]; \R^m)$ the space of the
$\R^m$-valued continuous processes $\{\varphi_t, t\in[0,T]\}$ such
that $\varphi_t$ is $\mathcal{H}_t$-measurable, $t\in[0,T]$, and
$\mathcal{M}^2(\mathbb{F}^W,[0,T];\R^m)$ the space of the
$\R^m$-valued square-integrable processes $\{\psi_t, t\in[0,T]\}$
such that $\psi_t$ is $\mathcal{F}_{[t,T]}^W$-measurable,
$t\in[0,T]$. Let $\HH_T^{\infty}(\R)$ be the set of
$\mathbb{H}$-progressively measurable processes which are almost
surely bounded by some real-valued $\F_T^B$-measurable random
variable, and let $\HH_T^2(\R^d)$ denote the set of all
$\R^d$-valued $\mathbb{H}$-progressively measurable processes
$\gamma=\{\gamma_t: t\in[0,T]\}$ such that
$\E\lt[\int^{T}_0|\gamma_t|^2\d t|\F^B_{T}\rt]<+\infty$, $\PP$-a.s.

\subsection{Russo-Vallois Integral}\label{sec:RV}
\setcounter{equation}{0} In a series of papers (\cite{Rv},
\cite{Rv1}, \cite{Rv2}, {\it etc.}), Russo and Vallois defined new
types of stochastic integrals, namely forward, backward and
symmetric integrals, which are extensions of the classical
Riemann-Stieltjes integral, and in fact these three integrals
coincide, when the integrator is a fBm with Hurst parameter
$H\in(1/2, 1)$. Here we will mainly use the backward Russo-Vallois
integral in this paper. It turns out to be a convenient definition
for stochastic integral with respect to our fBm $B$.

Let us recall some results by Russo and Vallois which we will use
later. In what follows, we make the convention that all continuous
processes $\{X_t, t\in[0,T]\}$ are extended to the whole line by
putting $X_t=X_0$, for $t<0$, and $X_t=X_T$, for $t>T$.

\begin{definition}{\label{def:rv}}
Let $X$ and $Y$ be two continuous processes. For $\ep>0$, we set
$$
I(\ep,t,X,\d Y)\triangleq\frac{1}{\ep}\int^t_0X(s)(Y(s)-Y(s-\ep))\d s,
$$
\[
C_\ep(X,Y)(t)\triangleq\frac{1}{\ep}\int^t_0(X(s)-X(s-\ep))(Y(s)-Y(s-\ep))\d
s, \ t\in[0,T].
\]
Then the backward Russo-Vallois integral is defined as the uniform
limit in probability as $\ep\to0^+$, if the limit exists. The
generalized bracket $[X,Y]$ is the uniform limit in probability of
$C_{\ep}(X,Y)$ as $\ep\to0^+$ $($of course, again under the
condition of existence$)$.
\end{definition}

We recall that (cf. Protter \cite{Pr}) a sequence of processes
$(H_n;n\ge0)$ converges to a process $H$ uniformly in probability if
\[
\lim_{n\to\infty}\mathbb{P}\lt(\sup_{t\in[0,T]}|H_n(t)-H(t)|>\alpha\rt)=0
\quad \mbox{for\ \ every}\ \ \alpha>0.
\]
In \cite{Rv2} (Theorem 2.1)  Russo and Vallois derived the It\^o
formula for the backward Russo-Vallois integral.
\begin{theorem}\label{thm:russo}
Let $f\in C^2(\R)$ and $X$ be a continuous process admitting the
generalized bracket, i.e., $[X,X]$ exists in the sense of Definition
\ref{def:rv}. Then for every $t\in[0,T]$, the backward Russo-Vallois
integral $\int^t_0f^\pr(X(s))\d X(s)$ exists and
\[
f(X(t))=f(X(0))+\int^t_0f^\pr(X(s))\d X(s)-\frac{1}{2}\int^t_0F^{\pr\pr}
(X(s))\d [X,X](s),
\]
for every $t\geq0$.
\end{theorem}

We list some properties of Russo-Vallois integral, which will be
used later in this paper.
\begin{proposition}\label{prop_rv1}
$(1)$. If $X$ is a finite quadratic variation process $($i.e.,
$[X,X]$ exists and $[X,X]_T<+\infty$, $\PP$-a.s.$)$ and $Y$ is a
zero quadratic variation process $($i.e., $[Y,Y]$ exists and equals
to zero$)$, then the mutual generalized bracket $[X,Y]$ exists and
vanishes, $\PP$-a.s.

$(2)$. If $X$ and $Y$ have $\PP$-a.s. H\"older continuous paths with
order $\alpha$ and $\beta$, respectively, such that $\alpha>0,$
$\beta>0$ and $\alpha+\beta>1$, then $[X,Y]=0$.

$(3)$. We assume that $X$ and $Y$ are continuous and admit a mutual
bracket. Then, for every continuous process $\{H(s):s\in[0,T]\}$,
$$
\int^\cdot_0H(s)\d C_\ep(X,Y)(s)\ \ \mbox{converges\  to}\ \
\int^\cdot_0H(s)\d [X,Y](s).
$$
\end{proposition}
The following proposition, which can be found in Russo and Vallois
\cite{Rv}, states the relationship between the Young integral (see
Young \cite{Yo}) and the backward Russo-Vallois integral.
\begin{proposition}
Let $X, Y$ be two real processes with paths being $\PP$-a.s. in
$C^\alpha$ and $C^\beta$, respectively, with $\alpha>0$, $\beta>0$
and $\alpha+\beta>1$. Then the backward Russo-Vallois integral
$\int^\cdot_0 Y\d X$ coincides with the Young integral
$\int^{\cdot}_0 Y\d^{(y)}X$.
\end{proposition}

\section{A Generalized It\^o Formula}
\setcounter{equation}{0} In this section we state a generalized
It\^o formula involving an It\^o backward integral with respect to
the Brownian motion $W$ and the Russo-Vallois integral with respect
to the fBm $B$. It will play an important role in our paper. It is
noteworthy that this It\^o formula corresponds to Lemma 1.3 in the
paper of Pardoux and Peng \cite{PP2} for the case of an fBm with
Hurst parameter $H=1/2$, i.e., when $B$ is a Brownian motion.

\begin{theorem}\label{thm:ito}
Let $\alpha\in C(\mathbb{H},[0,T]; \R)$ be a process of the form
\[
\alpha_t=\alpha_0+\int^t_0\beta_s\d s +\int^t_0 \gamma_s\dw\d W_s,
\quad t\in[0,T],
\]
where $\beta$ and $\gamma$ are $\mathbb{H}$-adapted processes and
$\PP\{\int^T_0|\beta_s|\d s<\infty\}=1$ and
$\PP\{\int^T_0|\gamma_s|^2\d s<+\infty\}=1$, respectively. Suppose
that $F\in C^2(\R\times\R)$. Then the Russo-Vallois integral
$\int_0^t\frac{\partial F}{\partial y}(\alpha_s,B_s)\d B_s $
$($defined as the uniform limit in probability of
$\frac{1}{\ep}\int^t_0(B_s-B_{s-\ep})\frac{\partial F}{\partial
y}(\alpha_s,B_s)\d s$$)$ exists for $0\le t\le T$,
and it holds that, $\mathbb{P}$-almost surely, for all $0\le t\le
T$,
\begin{equation}\label{eq:ito}
\begin{aligned}
F(\alpha_t,B_t)=&F(\alpha_0,0)+\int^t_0 \frac{\partial F}{\partial
x}(\alpha_s,B_s)\beta_s\d s +\int^t_0  \frac{\partial F}{\partial
x}(\alpha_s,B_s)\gamma_s\dw\d W_s + \int^t_0 \frac{\partial
F}{\partial
y}(\alpha_s,B_s)\d B_s\\
&-\frac12\int^t_0 \frac{\partial^2 F}{\partial
x^2}(\alpha_s,B_s)|\gamma_s|^2\d s.
\end{aligned}
\end{equation}
\end{theorem}
\noindent\textit{Proof:} {\bf{Step 1}}. First we suppose $F\in
C^2_b(\R\times \R)$ (i.e., the function $F$ is twice continuously
differentiable and has bounded derivatives of order less than or
equal to two) and there is a positive constant $C$ such that
$\int^T_0|\beta_s|\d s\le C$ and $\int^T_0|\gamma_s|^2\d s\le C$. It
is direct to check that
$$
F(\alpha_t,B_t)-F(\alpha_0,0)=\lim_{\ep\to0}\frac{1}{\ep}\int^t_0
(F(\alpha_s,B_s)-F(\alpha_{s-\ep},B_{s-\ep}))\d s.
$$
For simplicity we put $\alpha_{a,\ep,s}\triangleq
\alpha_{s}-a(\alpha_s-\alpha_{s-\ep})$ and $B_{a,\ep,s}\triangleq
B_{s}-a(B_s-B_{s-\ep})$, for any $a\in[0,1]$, $s\in[0,T],$ $\ep>0$.
We have
\begin{equation}\label{eq:3.2}
\begin{aligned}
&F(\alpha_s,B_s)-F(\alpha_{s-\ep},B_{s-\ep})\\
= &(\alpha_s-\alpha_{s-\ep}) \frac{\partial F}{\partial
x}(\alpha_s,B_s)+(B_s-B_{s-\ep}) \frac{\partial F}{\partial
y}(\alpha_s,B_s)-(\alpha_s-\alpha_{s-\ep})^2\int^1_0
\frac{\partial^2 F}{\partial x^2} (\alpha_{a,\ep,s},
B_{a,\ep,s})(1-a)\d
a\\
&-2(\alpha_s-\alpha_{s-\ep})(B_s-B_{s-\ep})\int^1_0 \frac{\partial^2
F}{\partial x\partial y} (\alpha_{a,\ep,s}, B_{a,\ep,s})(1-a)\d
a\\
&-(B_s-B_{s-\ep})^2\int^1_0 \frac{\partial^2 F}{\partial y^2}
(\alpha_{a,\ep,s}, B_{a,\ep,s})(1-a)\d a.
\end{aligned}
\end{equation}
By applying the stochastic Fubini theorem, we get that
$$
\begin{aligned}
&\frac{1}{\ep}\int^t_0 (\alpha_s-\alpha_{s-\ep})
\frac{\partial F}{\partial x}(\alpha_s,B_s)\d s\\
=&\frac{1}{\ep}\int^t_0 \lt(\int^s_{s-\ep}\beta_r\d r+\int^s_{s-\ep}
\gamma_r\dw\d W_r\rt)
\frac{\partial F}{\partial x}(\alpha_s,B_s)\d s\\
=&\frac{1}{\ep}\int^t_0 \int_r^{(r+\ep)\wedge
t}\beta_r\frac{\partial F}{\partial x}(\alpha_s,B_s)\d s\d r
+\frac{1}{\ep}\int^t_0 \int_r^{(r+\ep)\wedge t} \gamma_r
\frac{\partial F}{\partial x}(\alpha_s,B_s)\d s\dw\d W_r.
\end{aligned}
$$
Since $\frac{1}{\ep}\int_r^{(r+\ep)\wedge t}\frac{\partial
F}{\partial x}(\alpha_s,B_s)\d s$ is $\HH_t$-measurable and
converges to $\frac{\partial F}{\partial x}(\alpha_r,B_r)$ when
$\ep\to0$, it follows that, thanks to the continuity of
$\frac{\partial F}{\partial x}(\alpha_s,B_s)$,
$$
\begin{aligned}
&\lim_{\ep\to0} \sup_{t\in[0,T]}\lt|\int^t_0 \beta_r \lt(
\frac{1}{\ep}\int_r^{(r+\ep)\wedge t}\frac{\partial F}{\partial
x}(\alpha_s,B_s)\d s-\frac{\partial F}{\partial
x}(\alpha_r,B_r)\rt)\d r\rt|\\
\le& \lim_{\ep\to0} \sup_{r\in[0,T]}\lt|
\frac{1}{\ep}\int_r^{(r+\ep)\wedge t}\frac{\partial F}{\partial
x}(\alpha_s,B_s)\d s-\frac{\partial F}{\partial x}(\alpha_r,B_r)\rt|
\int^T_0 |\beta_r|\d r= 0,\ \ \mbox{in\ \ probability}.
\end{aligned}
$$
Thus, in virtue of the boundedness of $\frac{\partial F}{\partial
x}$, by the Dominated Convergence Theorem,
$$
\begin{aligned}
&\lim_{\ep\to
0}\mathbb{E}\lt[\lt.\sup_{t\in[0,T]}\lt|\int^t_0\gamma_r\lt(\frac{1}{\ep}\int_r^{(r+\ep)\wedge
t} \frac{\partial F}{\partial x}(\alpha_s,B_s)\d s-\frac{\partial
F}{\partial x}(\alpha_r,B_r)\rt)\dw\d W_r\rt|^2\rt|\F_T^B\rt]
\\
\le & \lim_{\ep\to
0}\mathbb{E}\lt[\lt.\sup_{r\in[0,T]}\lt|\frac{1}{\ep}\int_r^{(r+\ep)\wedge
t} \frac{\partial F}{\partial x}(\alpha_s,B_s)\d s-\frac{\partial
F}{\partial x}(\alpha_r,B_r)\rt|^2\int^t_0|\gamma_r|^2\d
r\rt|\F_T^B\rt]= 0, \ \PP-\textrm{a.s.}
\end{aligned}
$$
(Recall that $\lt(\frac{1}{\ep}\int_r^{(r+\ep)\wedge t}
\frac{\partial F}{\partial x}(\alpha_s,B_s)\d s-\frac{\partial
F}{\partial x}(\alpha_r,B_r)\rt)_{r\in[0,T]}$ is
$\mathbb{H}$-adapted and $W_\cdot-W_T$ is an $\mathbb{H}$-(backward)
Brownian motion.) Thus, we get
\begin{equation}\label{eq:3.3}
\begin{aligned}
\lim_{\ep\to0}\frac{1}{\ep}\int^t_0 (\alpha_s-\alpha_{s-\ep})
\frac{\partial F}{\partial x}(\alpha_s,B_s)\d s =\int^t_0\beta_r
\frac{\partial F}{\partial x}(\alpha_r,B_r)\d r +\int^t_0 \gamma_r
\frac{\partial F}{\partial x}(\alpha_r,B_r)\dw\d W_r,
\end{aligned}
\end{equation}
uniformly in probability. We notice that the generalized bracket of
$\alpha$ is the same as the classical one, i.e.,
$[\alpha,\alpha]_s=\int^s_0|\gamma_r|^2\d r$, $s\in[0,T]$.  We also
have
\begin{equation}\label{eq:3.4}
\begin{aligned}
&\frac{1}{\ep}\int^t_0 (\alpha_s-\alpha_{s-\ep})^2\int^1_0
\frac{\partial^2 F}{\partial x^2} (\alpha_{a,\ep,s},
B_{a,\ep,s})(1-a)\d a \d s\\
=&\frac{1}{2\ep}\int^t_0 (\alpha_s-\alpha_{s-\ep})^2
\frac{\partial^2 F}{\partial x^2} (\alpha_s, B_s)\d s+ A_{\ep,t}, \
\ t\in[0,T],
\end{aligned}
\end{equation}
where
$$
A_{\ep,t}=\frac{1}{\ep}\int^t_0 (\alpha_s-\alpha_{s-\ep})^2 \int^1_0
\lt(\frac{\partial^2 F}{\partial x^2} (\alpha_{a,\ep,s},
B_{a,\ep,s})-\frac{\partial^2 F}{\partial x^2} (\alpha_s,
B_s)\rt)(1-a)\d a \d s.
$$
Proposition \ref{prop_rv1} yields that $ \frac{1}{2\ep}\int^._0
(\alpha_s-\alpha_{s-\ep})^2 \frac{\partial^2 F}{\partial x^2}
(\alpha_s, B_s)\d s $ converges to $\int^._0 \frac{\partial^2
F}{\partial x^2} (\alpha_s, B_s)\d [\alpha,\alpha]_s$ and the
continuity of $\frac{\partial^2 F}{\partial x^2}$, $\alpha$ and $B$
implies that $A_{\ep,t}$ converges to zero. A similar argument shows
that
$$
\frac{1}{\ep}\int^._0
2(\alpha_s-\alpha_{s-\ep})(B_s-B_{s-\ep})\int^1_0 \frac{\partial^2
F}{\partial x\partial y} (\alpha_{a,\ep,s}, B_{a,\ep,s})(1-a)\d a \d
s$$ and the term
$$
\frac{1}{\ep}\int^._0(B_s-B_{s-\ep})^2\int^1_0 \frac{\partial^2
F}{\partial y^2} (\alpha_{a,\ep,s}, B_{a,\ep,s})(1-a)\d a\d s
$$
converge in probability, respectively, to
$$
2\int^._0\frac{\partial^2 F}{\partial x\partial y} (\alpha_s, B_s)\d
[\alpha,B]_s\ \ \mbox{and} \ \ \int^._0\frac{\partial^2 F}{\partial
x\partial y} (\alpha_s, B_s)\d [B,B]_s.
$$
However, these both latter expressions are zero due to the fact that
$H\in(1/2,1)$ and Proposition \ref{prop_rv1} (Observe that the fBm
has H\"older continuous paths of any positive order less than $H$
almost surely).

Combining the above results with (\ref{eq:3.2}), (\ref{eq:3.3}) and
(\ref{eq:3.4}), we see that
$$
\begin{aligned}
\lim_{\ep\to0}\frac{1}{\ep}\int^t_0(B_s-B_{s-\ep})\frac{\partial
F}{\partial y}(\alpha_s,B_s)\d
s=&F(\alpha_t,B_t)-F(\alpha_0,0)-\int^t_0 \frac{\partial F}{\partial
x}(\alpha_s,B_s)\beta_s\d s \\&-\int^t_0  \frac{\partial F}{\partial
x}(\alpha_s,B_s)\gamma_s\dw\d W_s  +\frac12\int^t_0 \frac{\partial^2
F}{\partial x^2}(\alpha_s,B_s)|\gamma_s|^2\d s,
\end{aligned}
$$
uniformly in $t\in[0,T]$, in probability. Consequently, the integral
$\int^t_0\frac{\partial F}{\partial y}(\alpha_s, B_s)\d B_s$ exists
in Russo-Vallois' sense (Recall Definition \ref{def:rv}) and we get
the It\^o formula (\ref{eq:ito}) for $F\in C^2_b(\R\times\R)$.

{\bf Step 2}. Now we deal with the general case that
$\PP\{\int^T_0|\beta_s|\d s<\infty \}=\PP\{\int^T_0|\gamma_s|^2\d
s<\infty\}=1$. For each $n\in\mathbb{N}$, we define a sequence of
$\mathbb{H}$-stopping times by $\tau_n=\sup\{t\le
T:\int^T_t|\beta_s|\d s+\int_t^T|\gamma_s|^2\d s>n\}\vee 0$, so we
know that the processes $\{\beta_t^n:=\beta_t
{\bf{1}}_{[\tau_n,T]}(t), t\in[0,T]\}$ and $\{\gamma^n_t:=\gamma_t
{\bf{1}}_{[\tau_n,T]}(t), t\in[0,T]\}$ satisfy
$\int^T_0|\beta^n_s|\d s\le n$ and $\int^T_0|\gamma^n_s|^2\d s\le
n,$ $\PP$-a.s. Furthermore, as $n\to\infty$, $\tau_n\to 0,$
$\PP$-a.s. We consider the It\^o formula for the process $
\alpha^n_t=\alpha_0+\int^t_0\beta^n_s\d s+ \int^t_0\gamma_s^n\dw \d
W_s, t\in[0,T]. $ Thanks to the result of {\bf{Step 1}}, we have
that, for every $n\in\mathbb{N}$,
$$
\begin{aligned}
F(\alpha_t^n,B_t)=&F(\alpha_0,0)+\int^t_0 \frac{\partial F}{\partial
x} (\alpha_s^n,B_s)\beta^n_s\d s +\int^t_0 \frac{\partial
F}{\partial y} (\alpha_s^n,B_s) \d B_s +\int^t_0 \frac{\partial
F}{\partial x}
(\alpha_s^n,B_s)\gamma_s\dw\d W_s\\
&-\frac12\int^t_0\frac{\partial^2 F}{\partial x^2}
(\alpha_s^n,B_s)|\gamma_s|^2\d s.
\end{aligned}
$$
Since $\alpha^n$ converges to $\alpha$ uniformly in probability on
$[0,T]$, by letting $n\to \infty$ in the above equation, we deduce
that $$\lim_{n\to\infty}\int^t_0\frac{\partial F}{\partial y}
(\alpha_s^n,B_s) \d B_s=\lim_{n\to\infty}\lim_{\ep\to
0}\int^t_0\frac{1}{\ep}\frac{\partial F}{\partial y}
(\alpha_s^n,B_s)(B_s-B_{s-\ep})\d s=\lim_{\ep\to
0}\int^t_0\frac{1}{\ep}\frac{\partial F}{\partial y}
(\alpha_s,B_s)(B_s-B_{s-\ep})\d s$$ exists, it is the Russo-Vallois
integral $\int^t_0\frac{\partial F}{\partial y} (\alpha_s,B_s)\d
B_s$ and it equals to
$$
\begin{aligned}
F(\alpha_t,B_t)-F(\alpha_0,0)-\int^t_0 \frac{\partial F}{\partial x}
(\alpha_s,B_s)\beta_s\d s -\int^t_0 \frac{\partial F}{\partial
x}(\alpha_s,B_s)\gamma_s\dw\d W_s +\frac12\int^t_0\frac{\partial^2
F}{\partial x^2} (\alpha_s,B_s)|\gamma_s|^2\d s.
\end{aligned}
$$

{\bf{Step 3}.} Finally we consider the case $F\in C^2(\R\times \R)$.
We let $\{\varphi_N\}_{N\in\N}$ be a sequence of infinitely
differentiable functions with compact support such that
$\varphi_N(x)=x$ for $\{(x_1,x_2):{\max}(|x_1|,|x_2|)\le N\},
N\in\N$. We set $F_N(x)=F(\varphi_N(x))$, so that $F_N(x)\in
C_b^2(\R\times \R)$ for every $N>0$. We notice that
$F_N(\alpha_\cdot,B_\cdot)$ and $F(\alpha_\cdot,B_\cdot)$ coincide
on the set $\Omega_N=\{ \omega\in\Omega: \sup_{0\le s\le
t}|\alpha_s|\le N,\sup_{0\le s\le t}|B_s|\le N\}$ and that
$\Omega=\bigcup_{N\ge1} \Omega_N$. Due to {\bf{Step 1}} and
{\bf{Step 2}}, for every $N$,  we have
\[
\begin{aligned}
F_N(\alpha_t,B_t)=&F_N(\alpha_0,0)+\int^t_0 \frac{\partial
F_N}{\partial x} (\alpha_s,B_s)\beta_s\d s +\int^t_0 \frac{\partial
F_N}{\partial y} (\alpha_s,B_s) \d B_s
+\int^t_0 \frac{\partial F_N}{\partial x} (\alpha_s,B_s)\gamma_s\dw\d W_s\\
&-\frac12\int^t_0\frac{\partial^2 F_N}{\partial x^2}
(\alpha_s,B_s)|\gamma_s|^2\d s, \ \ t\in[0,T].
\end{aligned}
\]
Therefore, for every $N$, it holds on $\Omega_N$  that
\[
\begin{aligned}
F(\alpha_t,B_t)=&F(\alpha_0,0)+\int^t_0 \frac{\partial F}{\partial
x} (\alpha_s,B_s)\beta_s\d s +\int^t_0 \frac{\partial F}{\partial y}
(\alpha_s,B_s) \d B_s +\int^t_0 \frac{\partial F}{\partial x}
(\alpha_s,B_s)\gamma_s\dw\d W_s\\
&-\frac12\int^t_0\frac{\partial^2 F}{\partial x^2}
(\alpha_s,B_s)|\gamma_s|^2\d s,\ \ t\in[0,T].
\end{aligned}
\]
Finally, by letting $N$ tend to $+\infty$, we get the wished result.
The proof is complete. \hfill$\cajita$

\section{Doss-Sussman Transformation of Fractional Backward Doubly Stochastic Differential Equations}
\setcounter{equation}{0} In what follows, we use the following
hypotheses:
\bigskip

(H1) The functions $\sigma: \R^n\to \R^{n\times d}$ and $b:
\R^{n}\to\R^n$ are Lipschitz continuous.
\bigskip

(H2) The function $f: [0,T]\times\R^n\times\R\times\R^d\to\R$ is
Lipschitz in $(x,y,z)\in\R^n\times\R\times\R^d\to\R$ with
$|f(t,0,0,0)|\le C$ uniformly in $t\in[0,T]$, the function $g:
\R^n\to \R $ belongs to $C^3_b(\R)$ and the function $\Phi$ is
bounded.
\bigskip

We fix an arbitrary $t\in[0,T]\subset\R^+$. Let $(X_s^{t,x})_{0\leq
s \leq t}$ be the unique solution of the following stochastic
differential equation:
\begin{equation}
\label{eq:x-sde} \left\{
\begin{array}{l}
\d X_s^{t,x}=- b(X_s^{t,x})\d
s -\sigma(X_s^{t,x})\dw\d  W_s , \quad s\in\left[0,t\right],\\
X_t^{t,x}=x.
\end{array}\right.
\end{equation}
Here the stochastic integral $\int^t_0\cdot\dw\d W_s$ is again
understood as the backward It\^o one. The condition (H1) guarantees
the existence and uniqueness of the solution $(X_s^{t,x})_{0\leq s
\leq t}$ in $\mathcal{M}^2(\mathbb{F}^W,[0,T];\R^n)$. Our aim is to
study the following backward doubly stochastic differential
equation:
\begin{equation}
\label{eq:x-bdsde}
U_s^{t,x}=\Phi(X_0^{t,x})+\int^s_0f(r,X_r^{t,x},
U_r^{t,x},V_r^{t,x})\d r + \int^s_0 g(U_r^{t,x})\d
B_r-\int^s_0V_r^{t,x}\dw\d W_r, \ s\in[0,t].
\end{equation}
We emphasiz that the integral with respect to the fBm $B$ is
interpreted in the Russo-Vallois sense, while the integral with
respect to the Brownian motion $W$ is the It\^o backward one. If $B$
is a standard Brownian motion, equation (\ref{eq:x-bdsde}) coincides
with the BDSDE which was first studied by Pardoux and Peng
\cite{PP2} in 1994 (apart of a time inversion).

Before we investigate the BDSDE (\ref{eq:x-bdsde}), we first give
the definition of its solution.
\begin{definition}
A solution of equation $(\ref{eq:x-bdsde})$ is a couple of processes
$(U_s^{t,x},V_s^{t,x})_{s\in[0,t]}$ such that:\\
1). $(U_s^{t,x},V_s^{t,x})_{s\in[0,t]}\in
\HH^{\infty}_t(\R)\times\HH_t^2(\R^d)$;\\
2). The Russo-Vallois integral $\int^\cdot_0g(U_r^{t,x})\d B_r$ is
well defined on $[0,t]$;\\
3). Equation $(\ref{eq:x-bdsde})$ holds $\PP$-a.s.

\end{definition}

Unlike the classical case, the lack of the semimartingale property
of the fBm $B$ gives an extra difficulty in solving BDSDE
$(\ref{eq:x-bdsde})$ directly. However, the work of Buckdahn and Ma
\cite{BM} indicates another possibility to investigate this
equation: by using the Doss-Sussman transformation. Let us develop
the idea. We denote by $\eta$ the stochastic flow which is the
unique solution of the following stochastic differential equation:
\begin{equation}
\label{eq:eta} \eta(t,y)=y+\int^t_0g(\eta(s,y))\d B_s, \ t\in[0,T],
\end{equation}
where the integral is interpreted in the sense of Russo-Vallois. The
solution of such a stochastic differential equation can be written
as $\eta(t,y)=\alpha(y,B_t)$ via Doss transformation (see, for
example, Z\"ahle \cite{Za2}), where $\alpha(y,z)$ is the solution of
the ordinary differential equation
\begin{equation}\label{eq:alpha}
\left\{
\begin{array}{l}
\frac{\partial\alpha}{\partial z} (y,z)=g(\alpha(y,z)), \ z\in\R,\\
\alpha(y,0)=y.
\end{array}\right.
\end{equation}
By the classical PDE theory we know that, for every $z\in\R$, the
mapping $y\mapsto\alpha(y,z)$ is a diffeomorphism over $\R$ and
$(y,z)\mapsto \alpha(y,z)$ is $C^2$. In particular, we can define
the $y$-inverse of $\alpha(y,z)$ and we denote it by $h(y,z)$, such
that we have $\alpha(h(y,z))=y$, $(y,z)\in\R\times\R$. Hence, it
follows that
$$
\frac{\partial\alpha}{\partial y}(h(y,z),z) \frac{\partial
h}{\partial y}(y,z)=1\,\ \mbox{and}\ \
\frac{\partial\alpha}{\partial z}(h(y,z),z)+
\frac{\partial\alpha}{\partial y}(h(y,z),z) \frac{\partial
h}{\partial z}(y,z)=0.
$$
Therefore,
\[\begin{aligned}
\frac{\partial h}{\partial
z}(y,z)&=-\left(\frac{\partial\alpha}{\partial
y}(h(y,z),z)\right)^{-1}\frac{\partial\alpha}{\partial z}(h(y,z),z)
=-\frac{\partial h}{\partial y}(y,z)\frac{\partial\alpha}{\partial
z}(h(y,z),z)=-\frac{\partial h}{\partial y}(y,z)g(y).
\end{aligned}
\]
As a direct consequence we have that also
$\eta(t,\cdot)=\alpha(\cdot,B_t):\R\to\R$ is a diffeomorphism and,
thus, we can define
$\mathcal{E}(t,y):=\eta(t,\cdot)^{-1}(y)=h(y,B_t),
(t,y)\in[0,T]\times\R$. Moreover, by the It\^o formula (Theorem
\ref{thm:ito}), we have
\[
\d\mathcal{E}(t,y)=\d h(y,B_t)=\frac{\partial }{\partial z} h(y,B_t)
\d B_t=-\frac{\partial}{\partial y}\mathcal{E}(t,y)g(y)\d B_t, \
t\in[0,T],
\]
i.e., the process $\mathcal{E}$ satisfies the following equation:
\begin{equation}
\label{eq:E} \mathcal{E}(t,y)=y-\int^t_0\frac{\partial}{\partial
y}\mathcal{E}(s,y)g(y) \d B_s, \ t\in[0,T].
\end{equation}
Furthermore, we have the following estimates for $\eta$ and
$\mathcal{E}$.

\begin{lemma}\label{prop_3.1}
There exists a constant $C>0$ depending only on the bound of $g$ and
its partial derivatives such that for $\xi=\eta, \mathcal{E}$, it
holds that,  $P$-a.s., for all $(t,y)$,
\[
\begin{array}{cc}
|\xi(t,y)|\le |y|+C|B_t|,
&\exp\{-C|B_t|\}\le\lt|\frac{\partial}{\partial y}\xi\rt|
\le \exp\{C|B_t|\},\\ \lt|\frac{\partial^2}{\partial y^2}\xi\rt|\le
\exp\{C|B_t|\}, &\left|\frac{\partial^3}{\partial y^3}\xi\right|\le
\exp\{C|B_t|\}.
\end{array}
\]
\end{lemma}

\noindent\textit{Proof:} The first three estimates are similar to
those in Buckdahn and Ma \cite{BM}. So we only prove the last one.
For this end we define $\gamma(\theta,y,z)=\alpha(y,\theta z)$, for
$(\theta,y,z)\in[0,1]\times\R\times\R^n$. It follows from
(\ref{eq:alpha}) that
$$
\gamma(\theta,y,z)=y+\int^{\theta}_0\langle
g(\gamma(r,y,z)),z\rangle\d r.
$$
By differentiating this latter equation, we obtain,
\begin{equation}\label{eq:alpha_yyy}
\left\{
\begin{array}{ll}
\frac{\partial^4\gamma}{\partial \theta\partial y^3}
(\theta,y,z)=&\frac{\partial^3 g }{\partial
y^3}(\gamma(\theta,y,z))z\left( \frac{\partial \gamma}{\partial
y}(\theta,y,z)\right)^3+3 \frac{\partial^2 g }{\partial
y^2}(\gamma(\theta,y,z))z\frac{\partial \gamma}{\partial
y}(\theta,y,z)\frac{\partial^2 \gamma}{\partial
y^2}(\theta,y,z)\\&+\frac{\partial g }{\partial
y}(\gamma(\theta,y,z))z
\frac{\partial^3 \gamma}{\partial y^3}(\theta,y,z);\\
\frac{\partial^3\gamma}{\partial y^3}(0,y,z)=0,
\end{array}\right.
\end{equation}
and from the variation of parameter formula it follows that
$$
\begin{aligned}
\frac{\partial^3}{\partial y^3}\gamma(1,y,z)=\int^1_0
\exp\left\{\int^1_u \frac{\partial g }{\partial y}(\gamma(v,y,z))z\d
v\right\}&\bigg(\frac{\partial^3 g }{\partial
y^3}(\gamma(u,y,z))z\left( \frac{\partial \gamma}{\partial
y}(u,y,z)\right)^3\\
& \quad+3 \frac{\partial^2 g }{\partial
y^2}(\gamma(u,y,z))z\frac{\partial \gamma}{\partial
y}(u,y,z)\frac{\partial^2 \gamma}{\partial y^2}(u,y,z)\bigg)\d u.
\end{aligned}
$$
Thus, by using the first three estimates of this Lemma, we get
\[
\left|\frac{\partial^3}{\partial
y^3}\eta(t,y)\right|=\left|\frac{\partial^3}{\partial
y^3}\alpha(y,B_t)\right|\le \exp\{C|B_t|\}.
\]
Hence we have completed the proof. \hfill$\cajita$

Lemma \ref{prop_3.1} plays an important role in the rest of the
paper thanks to the following lemma.

\begin{lemma}\label{lemma2}
For any $C\in\R$, we have
\[
\E\lt[\exp\lt\{C\sup_{s\in[0,T]}|B_s|\rt\}\rt]<\infty.
\]
\end{lemma}

\noindent\textit{Proof:} The proof is similar as Lemma 2.4 in
\cite{Jl} (and even easier), so we omit it. \hfill$\cajita$

We denote by $\wt{\Omega}^\pr$ the subspace of $\Omega^\pr$ such
that
$\wt{\Omega}^\pr=\lt\{\omega^\pr\in\Omega^\pr:\sup_{s\in[0,T]}|B_s|<\infty\rt\}$.
It is clear that $\PP^\pr(\wt{\Omega}^\pr)=1$.

We let $(Y^{t,x}, Z^{t,x})$ be the unique solution of the following
BSDE:
\begin{equation}
\label{eq:x-BSDE}
Y_s^{t,x}=\Phi(X_0^{t,x})+\int^s_0\tilde{f}(r,X_r^{t,x},Y_r^{t,x},Z_r^{t,x})\d
r -\int^s_0 Z_r^{t,x}\dw \d W_r,
\end{equation}
where
\[\begin{aligned}
\tilde{f}(t,x,y,z)=&\frac{1}{\frac{\partial}{\partial y}\eta
(t,y)}\left\{f\lt(t,x,\eta(t,y), \frac{\partial}{\partial
y}\eta(t,y)z\rt)
+\frac12\mathrm{tr}\lt[z^T\frac{\partial^2}{\partial
y^2}\eta(t,y)z\rt]\right\}.
\end{aligned}
\]
This BSDE, studied over $\Omega=\Omega^\pr\times\Omega^{\pr\pr}$ and
driven by the Brownian motion
$W_r(\omega)=W_r(\omega^{\pr\pr})=\omega^{\pr\pr}(r)$, $r\in[0,T]$,
can be interpreted as an $\omega^\pr$-pathwise BSDE, i.e., as a BSDE
over $\Omega^{\pr\pr}$, considered for every fixed
$\omega^\pr\in\Omega^\pr.$ However, subtleties of measurability make
us preferring to consider the BSDE over $\Omega$, with respect to
the filtration $\mathbb{H}$. We point out that the coefficient
$\tilde{f}$ has a quadratic growth in $z$, while the terminal value
is bounded. BSDEs of this type have been studied by Kobylanski
\cite{Ko}. We state an existence and uniqueness result for this kind
of BSDE, but with a slight adaptation to our framework. For this we
consider a driving coefficient $G$ satisfying the following
assumptions:

\bigskip
(H3) The coefficient $G: \Omega\times[0,T]\times \R\times\R^d$ is
measurable, for every fixed $(y,z)$, progressively measurable, with
respect to the backward filtration $\mathbb{H}$ and $G$ is
continuous in $(t,y,z)$;
\bigskip

(H4) There exists some real-valued $\F_T^B$-measurable random
variable $K:\Omega^\pr\to \R$ such that $|G(t,y,z)|\le K(1+|z|^2)$.
\bigskip

(H5) There exist real-valued $\F_T^B$-measurable random variables
$C>0$, $\ep>0$, and $\F_T^B\otimes\mathcal{B}([0,T])$-measurable
functions $k,$ $l_\ep: \Omega^\pr\times[0,T]\to \R$ such that
$$
\left|\frac{\partial G}{\partial z}(t,y,z)\right|\le k(t)+C|z|, \
\mathrm{for \ all} \ (t,y,z),\  \PP- \rm{a.s.,}
$$
$$
\frac{\partial G}{\partial y} (t,y,z)\le l_\ep(t)+\ep |z|^2, \
\mathrm{for\ all}\ (t,y,z),\  \PP-\rm{a.s.}
$$
\begin{remark}\label{rmk1}
Due to the Lemmata \ref{prop_3.1} and \ref{lemma2}, the function
$\tilde{f}$ in the equation $(\ref{eq:x-BSDE})$ satisfies
$(H3)-(H5)$. In particular, $K=\exp\{C\sup_{t\in[0,T]}|B_t|\}$ for
$C\in\R^+$ appropriately chosen.
\end{remark}

Adapting the results by Kobylanski \cite{Ko} (Theorem 2.3 and
Theorem 2.6), we can state the following:
\begin{theorem}\label{thm:ko}
Let $G$ be a driver such that hypotheses $(H3)$-$(H5)$ hold and let
$\xi$ be a real-valued $\mathcal{H}_0$-measurable random variable,
which is bounded by a real-valued $\F_T^B$-measurable random
variable. Then there exists a unique solution $(Y,Z)\in
\mathcal{H}_T^\infty(\R)\times\mathcal{H}^2_T(\R^d)$ of BSDE
\begin{equation}\label{eq:BSDE}
Y_t=\xi+\int^t_0G(s,Y_s,Z_s)\d s- \int^t_0 Z_s\dw \d W_s, \
t\in[0,T].
\end{equation}
Moreover, there exists a real-valued $\F_T^B$-measurable random
variable $C$ depending only on
$\esssup_{[0,T]\times\Omega^{\pr\pr}}$
$|Y_t(\omega^\pr,\omega^{\pr\pr})|$ and $K$, such that
\[
\mathbb{E}\lt[\int^T_0|Z_s|^2\d s\big|\F_T^B\rt]\le C,\
\PP-\rm{a.s.}
\]
\end{theorem}

\begin{remark}
The conditional expectation $\E\lt[\cdot|\F_T^B\rt]$ is here
understood in the generalized sense: if $\xi$ is a nonnegative
$\HH_0$-measurable random variable,
\[
\E[\xi|\F^B_T]:=\lim_{n\to\infty}\uparrow\E[\xi\wedge n|\F_T^B](\le
\infty)
\]
is a well defined $\F_T^B$-measurable random variable. If $\xi$ is
not nonnegative we decompose $\xi=\xi^+-\xi^-, \xi^+=\max\{\xi,0\},
\xi^-=-\min\{\xi,0\}$ and we put
$\E[\xi|\F_T^B]=\E[\xi^+|\F_T^B]-\E[\xi^-|\F_T^B]$ on
$\{\min\{\E[\xi^+|\F_T^B],\E[\xi^-|\F_T^B]\}<\infty\}$.
\end{remark}

\noindent\textit{Proof of Theorem \ref{thm:ko}:} We observe that the
Brownian motion $W$ possesses the (backward) martingale
representation with respect to the backward filtration $\mathbb{H}$,
i.e., given an $\HH_0$-measurable random variable $\xi$ such that
$\E\lt[\xi^2|\F^B_T\rt]<\infty$, $\PP$-a.s., there exists a unique
process $\gamma\in \HH_T^2(\R^d)$ such that
\[
\xi=\E[\xi|\F_T^B]-\int^T_0\gamma_r\dw \d W_r, \ \PP-\rm{a.s.}
\]
This martingale representation property allows to show the existence
and uniqueness of a solution of (\ref{eq:BSDE}) when $G$ is of
linear growth and Lipschitz in $(y,z)$. Combining this with the
approach by Kobylanski \cite{Ko} allows to obtain the result stated
in Theorem \ref{thm:ko}.\hfill$\cajita$

\bigskip

By using Theorem $\ref{thm:ko}$, we are now able to characterize
more precisely the solution of BSDE (\ref{eq:x-BSDE}).

\begin{theorem}\label{thm:sol}
Under our standard assumptions on the coefficients $\sigma, b , f$
and $\Phi$, BSDE $(\ref{eq:x-BSDE})$ admits a unique solution
$(Y^{t,x}, Z^{t,x})$ in $\HH_t^\infty(\R)\times\HH_t^2(\R^d)$.
Moreover, there exists a positive increasing process $\theta\in
L^0(\mathbb{H},\R)$ such that $|Y^{t,x}_s|\le \theta_s$,
$\E\lt[\int^\tau_0|Z^{t,x}_s|^2\d s |\HH_\tau\rt]\le
\exp\{\exp\{C\sup_{s\in[0,t]}|B_s|\}\}$, $\PP$-a.s., for all
$\mathbb{H}$-stopping times $\tau$ $(0\le \tau\le t)$, where $C$ is
a constant chosen in an adequate way. Furthermore, the process
$(Y^{t,x},Z^{t,x})$ is $\mathbb{G}$-adapted.
\end{theorem}

\noindent\textit{Proof:} Due to Theorem \ref{thm:ko}, equation
(\ref{eq:x-BSDE}) has a unique solution $(Y^{t,x}, Z^{t,x})$ in
$\HH_t^\infty(\R)\times\HH_t^2(\R^d)$.

{\bf Step 1}. In order to give the estimates of $(Y^{t,x},
Z^{t,x})$, we proceed as in Lemma 5.3 in \cite{BM}. In particular,
we can show that there exists an increasing positive process
$\theta\in L^0(\mathbb{F}^B, [0,T])$ such that $\PP$-a.s.,
$|Y^{t,x}_s|\le \theta_s,  0\le s \le t \le T$. This process
$(\theta_s)$ can be chosen as the solution of the following ordinary
differential equation
\[
\frac{\d \theta_s}{\d s}=\exp\{C\sup_{0\le r\le
t}|B_r|\}(1+\theta_s);\ \theta(0)=|\Phi(X_0^{t,x})|,
\]
for some suitably chosen real constant $C$, i.e.,
$$
\theta_s=(|\Phi(X_0^{t,x})|+1)\exp\lt\{\exp\lt\{C\sup_{0\le r\le
t}|B_r|\rt\}s\rt\}-1.
$$

Indeed, for $M>0$, let $\varphi_M(y)$ be a $C^\infty$ function such
that $0\le \varphi_M\le 1$, $\varphi_M(y)=1$ for $|y|\le M$ and
$\varphi_M(y)=0$ for $|y|\ge M+1$. Defining a new function
$\tilde{f}^M$ by $\tilde{f}^M(t,x,y,z)\triangleq
\tilde{f}(t,x,y,z)\varphi_M(y)$ we see that the function
$\tilde{f}^M$ also satisfies conditions (H3)-(H5). According to
Theorem \ref{thm:ko}, there exists a unique solution
$(Y^{M,t,x},Z^{M,t,x})$ of equation (\ref{eq:x-BSDE}) with
$\tilde{f}$ being replaced with $\tilde{f}^M$. The stability result
in Kobylanski \cite{Ko} shows that, when $M\to +\infty$, there
exists a subsequence of $Y^{M,t,x}$ converging to $Y^{t,x}$
uniformly in probability. Therefore, following Buckdahn and Ma's
approach and slightly adapted, we only need to prove that
$Y^{M,t,x}_s$ is uniformly bounded by $\theta_s$. We apply the
Tanaka formula to $|Y^{M,t,x}|$ to get that
\[
\lt|Y^{M,t,x}_s\rt|=|\Phi(X_0^{t,x})|+\int^s_0\mathrm{sgn}\lt(Y^{M,t,x}_r\rt)
\tilde{f}^M(r,X_r^{t,x},Y^{M,t,x}_r,Z^{M,t,x}_r)\d r-
\int^s_0\mathrm{sgn}\lt(Y^{M,t,x}_r\rt)Z^{M,t,x}_r\dw \d W_r
+L_s-L_0,
\]
for a local-time-like process $L$ such that $L_t=0$ and
$L_s=\int^t_s1_{\{Y^{M,t,x}_r=0\}}\d L_r$. Now we define a new
function $\psi(y)$ by letting $\psi(y)=e^{2Ky}-1-2Ky-2K^2y^2$, for
$y>0$, and $\psi(y)=0$ for $y\le 0$, where $K$ is the bound in
Remark \ref{rmk1}. Then we apply the It\^o formula to
$\psi(\lt|Y^{M,t,x}_s\rt|-\theta_s)$ and get
\begin{equation}\label{eq:YM}
\begin{aligned}
&\psi\lt(\lt|Y^{M,t,x}_s\rt|-\theta_s\rt)\\=&
\int^s_0\psi^\pr\lt(\lt|Y^{M,t,x}_r\rt|-\theta_r\rt)
\mathrm{sgn}\lt(Y^{M,t,x}_r\rt)
\lt(\tilde{f}^M\lt(r,X_r^{t,x},Y^{M,t,x}_r,Z^{M,t,x}_r\rt)
-\exp\{C\sup_{0\le u\le t}|B_u|\}(1+\theta_r)\rt)\d
r\\&-\int^s_0\psi^\pr\lt(\lt|Y^{M,t,x}_r\rt|-\theta_r\rt)
\mathrm{sgn}\lt(Y^{M,t,x}_r\rt)Z^{M,t,x}_r\dw \d W_r
+\int^s_0\psi^\pr\lt(\lt|Y^{M,t,x}_r\rt|-\theta_r\rt)\d L_r\\&
-\frac12\int^s_0\psi^{\pr\pr}\lt(\lt|Y^{M,t,x}_r\rt|-\theta_r\rt)\lt|Z^{M,t,x}_r\rt|^2\d
r
\end{aligned}
\end{equation}
The property of $\psi$ shows that
$\int^s_0\psi^\pr\lt(\lt|Y^{M,t,x}_r\rt|-\theta_r\rt)\d
L_r=\int^s_0\psi^\pr\lt(-\theta_r\rt)\d L_r=0$. We also have
\[
\begin{aligned}
&\int^s_0\psi^\pr\lt(\lt|Y^{M,t,x}_r\rt|-\theta_r\rt)
\mathrm{sgn}\lt(Y^{M,t,x}_r\rt)
\lt(\tilde{f}^M\lt(r,X_r^{t,x},Y^{M,t,x}_r,Z^{M,t,x}_r\rt)
-\exp\{C\sup_{0\le u\le t}|B_u|\}(1+\theta_r)\rt)\d r\\
\le &\int^s_0 \psi^\pr\lt(\lt|Y^{M,t,x}_r\rt|-\theta_r\rt)
\lt(K\lt(\varphi(Y^{M,t,x}_r)\lt|Y^{M,t,x}_r\rt|-\theta_r\rt)+K|Z^{M,t,x}_r|^2\rt)\d
r.
\end{aligned}
\]
Since $\psi^{\pr\pr}-2K\psi^\pr\ge0$, we get from equation
(\ref{eq:YM}) that
\[
\begin{aligned}
\psi\lt(\lt|Y^{M,t,x}_s\rt|-\theta_s\rt)\le& K
\int^s_0\psi^\pr\lt(\lt|Y^{M,t,x}_r\rt|-\theta_r\rt)
\lt(\varphi(Y^{M,t,x}_r)\lt|Y^{M,t,x}_r\rt|-\theta_r\rt)\d
r\\&-\int^s_0\psi^\pr\lt(\lt|Y^{M,t,x}_r\rt|-\theta_r\rt)
\mathrm{sgn}\lt(Y^{M,t,x}_r\rt)Z^{M,t,x}_r\dw \d W_r.
\end{aligned}
\]
Thus, we deduce that
\[
\begin{aligned}
\E\lt[\lt.\psi\lt(\lt|Y^{M,t,x}_s\rt|-\theta_s\rt)\rt|\HH_s\rt]\le
&\E\lt[\lt.\ K \int^s_0\psi^\pr\lt(\lt|Y^{M,t,x}_r\rt|-\theta_r\rt)
\lt(\varphi(Y^{M,t,x}_r)\lt|Y^{M,t,x}_r\rt|-\theta_r\rt)\d
r\rt|\HH_s\rt].
\end{aligned}
\]
From the definition of $\psi$ we get that
$\psi^\pr(y)=2K(\psi(y)+2K^2y^2)$. Hence, we have
\[
\begin{aligned}
&\E\lt[\lt.\psi\lt(\lt|Y^{M,t,x}_s\rt|-\theta_s\rt)\rt|\HH_s\rt]\\
\le&\E\lt[\lt.\ \int^s_02K^2\psi\lt(\lt|Y^{M,t,x}_r\rt|-\theta_r\rt)
\lt(\varphi(Y^{M,t,x}_r)\lt|Y^{M,t,x}_r\rt|-\theta_r\rt)+4K^4
\lt(\lt|Y^{M,t,x}_r\rt|-\theta_r\rt)^3\d r\rt|\HH_s\rt].
\end{aligned}
\]
There also exists a $\tilde{K}$ such that $y^3\le\tilde{K}\psi(y)$,
for all $y\in\R$. Consequently, we get that
\[
\begin{aligned}
\E\lt[\lt.\psi\lt(\lt|Y^{M,t,x}_s\rt|-\theta_s\rt)\rt|\HH_s\rt]
\le2K^2(M+\|\theta\|_{\infty,[0,t]}+2K^2\tilde{K})\int^s_0
\E\lt[\lt.\psi\lt(\lt|Y^{M,t,x}_r\rt|-\theta_r\rt)\rt|\HH_s\rt]\d r.
\end{aligned}
\]
Finally, the Gronwall inequality shows that
$\psi\lt(\lt|Y^{M,t,x}_s\rt|-\theta_s\rt)=0$,  for any $s\in[0,t]$,
$\PP$-a.s. Therefore, $\lt|Y^{M,t,x}_s\rt|\le\theta_s$, for any
$s\in[0,t]$, $\PP$-a.s.

{\bf Step 2}. We apply the It\^o formula to $e^{aY^{t,x}_s}$, with
$a$ being a real-valued $\F_T^B$-measurable random variable to be
determined later, and we obtain
\[
\begin{aligned}
e^{aY^{t,x}_s}&=e^{a\Phi(X_0^{t,x})}+\int^s_0ae^{aY_s}\tilde{f}(r,X^{t,x}_r,Y^{t,x}_r,
Z^{t,x}_r)\d r-\int^s_0\frac{1}{2}a^2e^{aY^{t,x}_r}|Z^{t,x}_r|^2\d
r-\int^s_0ae^{aY^{t,x}_r}Z^{t,x}_r\dw
\d W_r\\
&\le e^{a\Phi(X_0^{t,x})}+\int^s_0\left(-\frac{1}{2}a^2+|a|K\right)
e^{aY^{t,x}_r}|Y^{t,x}_r|^2\d r+\int^s_0 |a|K e^{aY^{t,x}_r}\d
r-\int^s_0ae^{aY^{t,x}_r}Z^{t,x}_r\dw \d W_r.
\end{aligned}
\]
Let
$\zeta(\omega^\pr):=\esssup_{[0,t]\times\Omega^{\pr\pr}}|Y^{t,x}_s
(\omega^\pr,\omega^{\pr\pr})|(<+\infty,
\omega^\pr\in\wt{\Omega}^\pr)$. Hence, by taking the conditional
expectations $\E[\cdot|\HH_t]$ on both sides, we deduce that for any
$\mathbb{H}$-stopping times $\tau\in[0,T]$,
\[
\begin{aligned}
&\left(\frac12a^2-|a|K\right)e^{-|a|\zeta}
\mathbb{E}\left[\int^\tau_0 |Z^{t,x}_r|^2\d
r\bigg|\HH_\tau\right]\le \left(\frac12a^2-|a|K\right)
\mathbb{E}\left[\left.\int^\tau_0 e^{aY^{t,x}_r}|Z_r|^2\d r
\rt|\HH_\tau\right]
\\
\le&\mathbb{E}\left[\left.e^{a\Phi(X_0^{t,x})}-e^{aY^{t,x}_\tau}+\int^\tau_0
|a|K e^{aY^{t,x}_r}\d s\right|\HH_\tau\right].
\end{aligned}
\]
We can choose $a=4K$ such that $\frac12a^2-|a|K=4K^2$ and we get,
keeping in mind that here $K$ is a random variable bounded by
$\exp\{C\sup_{s\in[0,t]}|B_s|\}$ ($C\in\R^+$ is a real constant, see
Remark \ref{rmk1}),
\begin{equation}\label{est-Z}
\mathbb{E}\left[\left.\int^\tau_0 |Z^{t,x}_r|^2\d r\right|
\HH_\tau\right]\le \frac{(2+4K^2t)}{4K^2}e^{8K}\le
\exp\lt\{\exp\lt\{C\sup_{s\in[0,t]}|B_s|\rt\}\rt\}.
\end{equation}

{\bf Step 3}. Let us show that the process $(Y^{t,x}, Z^{t,x})$ is
not only $\mathbb{H}$- but also $\mathbb{G}$-adapted. For this we
consider for an arbitrarily given $\tau\in[0,t]$ equation
(\ref{eq:x-BSDE}) over the time interval $[0,\tau]$:
\begin{equation}\label{x-bsde2}
Y_s^{t,x}= \Phi(X_0^{t,x}) +\int^s_0 \tilde{f}\left(r, X_r^{t,x},
Y_r^{t,x}, Z_r^{t,x}\right)\d r -\int_0^sZ_r^{t,x}\dw \d W_r, \quad
s\in\left[0,\tau\right].
\end{equation}
Let
$\mathcal{H}_t^{\tau}:=\mathcal{F}_{t,T}^W\vee\mathcal{F}_{\tau}^{B},$
$t\in[0,\tau]$. Then
$\mathbb{H}^{\tau}=\{\mathcal{H}_t^{\tau}\}_{t\in[0,\tau]}$ is a
backward Brownian filtration enlarged by a $\sigma$-algebra
generated by the fBm $B$, which is independent of the Brownian
filtration. Thus, with respect to $\mathbb{H}^\tau$, the Brownian
motion $W$ has the martingale representation property. Since
$\tilde{f}(r,x,y,z)$ is $\mathcal{G}_t$- and, hence, also
$\mathcal{H}_t^{\tau}$-measurable, $\d r$ a.e. on $[0,\tau]$, it
follows from the classical BSDE theory (or Theorem \ref{thm:ko})
that BSDE (\ref{x-bsde2}) admits a unique solution
$\left(Y^{t,x,\tau},Z^{t,x,\tau}\right) \in
\HH_\tau^\infty(\R)\times\HH_\tau^2(\R^d)$. On the other hand, also
$\left(Y^{t,x,\tau}_r,Z^{t,x,\tau}_r\right)_{r\in[0,\tau]}$ is a
solution of (\ref{eq:x-BSDE}). Hence,
$\left(Y^{t,x}_r,Z^{t,x}_r\right)=
\left(Y^{t,x,\tau}_r,Z^{t,x,\tau}_r\right)$, $\d r$ a.e., for
$t<\tau$. Consequently, $\left(Y^{t,x}_r,Z^{t,x}_r\right)$ is
$\mathcal{H}_r^\tau$- measurable, $\d r$ a.e., for $r<\tau.$
Therefore, letting $\tau\dw t$ we can deduce from the right
continuity of the filtration $\mathbb{F}^B$ that
$\left(Y^{t,x},Z^{t,x}\right)$ is
$\mathbb{G}$-adapted.\hfill$\cajita$

\begin{remark}
We remind the reader that the bound we get in $(\ref{est-Z})$ is
only $\PP$-a.s. finite, but not square-integrable. As a matter of
fact, it is hard to prove directly that $Z^{t,x}$ is a
square-integrable process, which constitutes the main reason that we
use instead the space $\HH_t^2(\R^d)$. Hence, the major difference
between our work and Buckdahn and Ma \cite{BM} is: In the classical
case, a priori we can solve the BDSDE in the first step to get the
square integrability of $Z$, but in the fractional case, there is
not a direct way to solve the BDSDE.
\end{remark}

Now we are ready to give the main result of this section by linking
the BSDEs (\ref{eq:x-bdsde}) and (\ref{eq:x-BSDE}) with the help of
the Doss-Sussman transformation.
\begin{theorem}\label{thm:ds}
Let us define a new pair of processes $(U^{t,x}, V^{t,x})$ by
$$
U_s^{t,x}=\eta(s, Y_s^{t,x}), \ \ V_s^{t,x}=\frac{\partial}{\partial
y}\eta(s, Y_s^{t,x})Z_s^{t,x},
$$
where $(Y^{t,x},Z^{t,x})$ is the solution of BSDE
$(\ref{eq:x-BSDE})$. Then $(U^{t,x}, V^{t,x})\in
\HH_t^{\infty}(\R)\times\HH_t^2(\R^d)$ is the solution of BDSDE
$(\ref{eq:x-bdsde})$.
\end{theorem}
\begin{remark}
The above theorem can be considered as a counterpart of Theorem 3.9
in Jing and Le\'on \cite{Jl} for the semi-linear case when $H<1/2$.
However, since we use here a different Hurst parameter $H$ and a
different type of stochastic integral with respect to fBm $B$, the
above theorem obviously does not cover the result in \cite{Jl}.
\end{remark}

\noindent\textit{Proof of Theorem \ref{thm:ds}:} The fact that
$(U^{t,x}, V^{t,x})\in \HH_t^{\infty}(\R)\times\HH_t^2(\R^d)$
follows directly from Theorem \ref{thm:sol} and Lemma
\ref{prop_3.1}. In order to prove the remaining part of the theorem,
we just have to apply the It\^o formula (Theorem \ref{thm:ito}) to
$\alpha(Y_s^{t,x},B_s)$, noticing $\alpha(Y_s^{t,x},B_s)=\eta(s,
Y_s^{t,x})$, to obtain that, for $(s,x)\in[0,t]\times\R^d$, the
Russso-Vallois integral $\int^s_0g(U^{t,x}_r)\d B_r$ exists and
\begin{equation}
\begin{aligned}
U_s^{t,x}=&\Phi_0(X_0^{t,x}) +\int^s_0 \frac{\partial}{\partial
y}\alpha( Y_r^{t,x},B_r)\bigg[\lt({\frac{\partial}{\partial y}\eta
(r, Y_r^{t,x})}\rt)^{-1} \bigg\{f\left(r,X_r^{t,x},\eta(s,
Y_r^{t,x}), \frac{\partial}{\partial y}\eta(r,
Y_r^{t,x})Z_r^{t,x}\right)
\\
& +\frac12\mathrm{tr}\lt[(Z^{t,x}_r)^T\frac{\partial^2}{\partial
y^2}\eta(r, Y_r^{t,x})Z_r^{t,x}\rt]\bigg\}\d r-Z_r^{t,x}\dw\d
W_r\bigg]
-\frac12\int^s_0\mathrm{tr}\lt[(Z^{t,x}_r)^T\frac{\partial^2}{\partial
y^2}\eta(r, Y_r^{t,x})Z_r^{t,x}\rt]\d
r\\
&+\int^s_0 \frac{\partial}{\partial z}\alpha( Y_r^{t,x},B_r) \d B_r.
\end{aligned}
\end{equation}
Consequently,
\begin{equation}
\begin{aligned}
U_s^{t,x}=\Phi(X_0^{t,x})+\int^s_0
f(r,X_r^{t,x},U_r^{t,x},V_r^{t,x})\d r + \int^s_0 g(U_r^{t,x}) \d
B_r-\int^s_0V_r^{t,x}\dw \d W_r.
\end{aligned}
\end{equation}
The proof is complete.\hfill$\cajita$

Now we close this section by a simple example to illustrate the idea
of the above procedure.
\begin{example}
Let us consider the linear case $($ for simplicity of notations we
omit the superscript $(t,x)$$)$:
\begin{equation}\label{eq:example}
\left\{
\begin{array}{l}
\d U_s= U_s\d B_s+f(U_s,V_s)\d s +V_s\dw \d W_s, \ s\in[0,t];\\
U_0=\Phi(X_0^{t,x}).
\end{array} \right.
\end{equation}
It is elementary to show that $\alpha(y,z)=ye^z$. Hence, the
solution $(U,V)$ of equation $(\ref{eq:example})$ is given by
$U_s=Y_se^{B_s}$ and $V_s=Z_se^{B_s}$, where $(Y,Z)$ is the solution
of
\begin{equation}\label{eq:example1}
\left\{
\begin{array}{l}
\d Y_s= \hat{f}(Y_s,Z_s)\d s +Z_s\dw \d W_s, \ s\in[0,t];\\
Y_0=\Phi(X_0^{t,x}).
\end{array} \right.
\end{equation}
and the function $\hat{f}$ is defined by
$\hat{f}(y,x)=e^{-B_t}f(ye^{B_t},ze^{B_t})$. Obviously $\hat{f}$ is
Lipschitz in $(y,z)$ as far as $f$ is Lipschitz. By following a
classical Malliavin calculus method $($see, e.g., Pardoux and Peng
\cite{PP2}$)$, we can get that $Z$ is $\PP$-a.s. uniformly bounded
and, thus, the process $\{Y_s,s\in[0,t]\}$ is $(1/2-\ep)$-H\"older
continuous in $s$, for all $\ep\in (0,1/2)$. Therefore, we have, for
every $r,s\in [0,t]$,
\[
\begin{aligned}
&|U_r-U_s|=|Y_re^{B_r}-Y_se^{B_s}|\le
\|Y\|_{\infty}|e^{B_r}-e^{B_s}|+ e^{|B_s|}|Y_r-Y_s|\\
\le&\exp\lt\{\sup_{0\le u\le
t}|B_u|\rt\}\|Y\|_{\infty}C(\omega^\pr)|r-s|^\alpha+\exp\lt\{\sup_{0\le
u\le t}|B_u|\rt\}C(\omega^\pr)|r-s|^{1/2-\ep}, \rm{for}\ \rm{a.a.}
\omega^\pr\in \Omega^\pr,
\end{aligned}
\]
where we can choose $\alpha$ to be in $(1/2,H)$. That is to say, we
can choose $0< \ep< \alpha-1/2$ and get that  the process $\{U_s,
s\in [0,t]\}$ has $\alpha$-H\"older continuous paths. Hence, instead
of using the Russo-Vallois integral with respect to the fractional
Brownian motion in equation $(\ref{eq:example})$ we can use the
classical Young integral. Furthermore, thanks to Proposition
\ref{prop_rv1}, these two integrals coincide.
\end{example}

\begin{remark}\label{rmk2}
\em{In the latter example, the Young integral and the Russo-Vallois
integral coincide. Unfortunately, the H\"older continuity seems to
be very hard to deduce in the general, nonlinear case. As a
consequence, we have to work with the more general Russo-Vallois
integral.}
\end{remark}

\section{Associated Stochastic Partial Differential Equations}
\setcounter{equation}{0}

In this section we discuss briefly the relationship between an
associated SPDE and a PDE with stochastic coeffients. For
simplicity, we only show the relationship for the case of classical
solutions. For a complete discussion of the case of (stochastic)
viscosity solutions, one can proceed by adapting the approaches in
Buckdahn and Ma \cite{BM} and Jing and Le\'on \cite{Jl}.

Let $\mathscr{L}$ be the second order elliptic differential
operator:
\[
\mathscr{L}=\frac12\sum_{i,j=1}^n (\sigma\sigma^T)_{ij}(x)
\frac{\partial^2 }{\partial x_i \partial x_j}+ \sum_{i=1}^n b_i(x)
\frac{\partial }{\partial x_i},
\]
which means that it is the infinitesimal generator of the Markovian
process $\{X^{t,x}_s, s\in[0,t]\}$ defined by equation
(\ref{eq:x-sde}).

Our aim is to study the following semilinear SPDE driven by the fBm
$B$:
\begin{equation}\label{spde}
\left\{
\begin{array}{lc}
\d v(t,x) = \left[\mathscr{L} v (t,x) -
{f}(t,x,v(t,x),\sigma(x)^T\frac{\partial}{\partial x}v(t,x))\rt]\d t
+ g(v(t,x))\d B_t,& (t,x)\in (0,T)\times
\R^n;\\
v(0,x)=\Phi(x),&x\in\R^n.
\end{array}\right.
\end{equation}
In the case that $B$ is a Brownian motion (an fBm with Hurst
parameter $H=1/2$), it is well known (see, Pardoux and Peng
\cite{PP2}, Buckdahn and Ma \cite{BM}) that the random field
$v(t,x):=U^{t,x}_t$ solves $(\ref{spde})$ (in the viscosity sense if
the coefficients are Lipschitz, and in the classical sense if the
coefficients are $C^3_b$), where $U^{t,x}$ is the solution of BDSDE
(\ref{eq:x-bdsde}). Thus, it is natural to raise the following
question: Can we also solve the SPDE (\ref{spde}) driven by the fBm
$B$, by studying the properties of the solution of its associated
BDSDE (\ref{eq:x-bdsde})? The answer is positive. Indeed, by
applying the Doss-Sussman transformation, we will show that the PDE:
\begin{equation}\label{w-pde}
\left\{
\begin{array}{cllc}
\d u(t,x)&=&\lt(\mathscr{L} u (t,x) -
\tilde{f}(t,x,u(t,x),\sigma(x)^T\frac{\partial}{\partial
x}u(t,x))\rt)\d t,& (t,x)\in (0,T)\times
\R^n;\\
u(0,x)&=&\Phi(x),&x\in\R^n,
\end{array}\right.
\end{equation}
is transformed into SPDE (\ref{spde}), where we recall that
$$
\tilde{f}(t,x,y,z)={\lt(\frac{\partial}{\partial y} \eta
(t,y)\rt)^{-1}}\left\{f(t,x,\eta(t,y), \frac{\partial}{\partial
y}\eta(t,y)z) +\frac12\mathrm{Tr}\lt[z\frac{\partial^2}{\partial
y^2}\eta(t,y)z^T\rt]\right\}.
$$
We also observe that following a similar argument as in Kobylanski
\cite{Ko}, under some smoothness assumptions, $u(t,x):=Y^{t,x}_t,
(t,x)\in[0,T]\times\R^n$ is the solution of equation (\ref{w-pde}),
where $Y^{t,x}$ is the solution of BSDE (\ref{eq:x-BSDE}).

First we give the definition for the classical solutions of
equations (\ref{spde}) and (\ref{w-pde}).
\begin{definition}
We say a stochastic field $w: \Omega^\pr\times[0,T]\times\R^n\to\R$
is a classical solution of equation $(\ref{spde})$ $($resp.,
$(\ref{w-pde}))$, if $w\in C^{0,2}_{\mathbb{F}^B}$
\footnote{$C^{0,2}_{\mathbb{F}^B}$ is the space of
${\mathbb{F}^B}$-adapted processes $u(t,x)$ which are continuous in
$t$ and $C^2$ in $x$.} and satisfies equation $(\ref{spde})$
$($resp., $(\ref{w-pde}))$.
\end{definition}

We have the following proposition.
\begin{proposition}\label{prop_u}
Suppose that $u$ is a classical solution of equation
$(\ref{w-pde})$, then
$\hat{u}(t,x)\triangleq\eta(t,u(t,x))=\alpha(u(t,x),B_t)$ is  a
classical solution of SPDE $(\ref{spde})$. The converse holds also
true: Every classical solution $\hat{u}$ of equation $(\ref{spde})$
defines a classical solution $u(t,x)=\mathcal{E}(t,\hat{u}(t,x))$ of
equation $(\ref{w-pde})$.
\end{proposition}
\noindent\textit{Proof:} The claim that $\hat{u}(t,x)\in
C^{0,2}_{\mathbb{F}^B}$ follows from the regularity property of the
functions $\alpha$ and $u$ and the fact that $u$ is
$\mathbb{F}^B$-adapted. Moreover, we first observe that
$\hat{u}(0,x)=\eta(0,u(0,x))=u(0,x)=\Phi(x)$, and we apply the It\^o
formula to $\hat{u}(t,x)$ to obtain:
\begin{equation}\label{eq:hatu}
\begin{aligned}
\d \hat{u}(t,x) =& \frac{\partial}{\partial
y}\alpha(u(t,x),B_t)\lt(\mathscr{L} u (t,x) -
\tilde{f}\lt(t,x,u(t,x),\sigma(x)^T\frac{\partial}{\partial
x}u(t,x)\rt)\rt)\d t
+\frac{\partial}{\partial z}\alpha(u(t,x),B_t)\d B_t\\
=& \Bigg[\frac{\partial}{\partial y}\alpha(u(t,x),B_t)\mathscr{L} u
(t,x)- f\lt(t,x,\eta(t,u(t,x)),\frac{\partial}{\partial
y}\eta(t,u(t,x))\sigma(x)^T\frac{\partial}{\partial x}u(t,x)\rt)\\
&\quad-\frac12\mathrm{Tr}\lt[\sigma(x)^T\frac{\partial}{\partial
x}u(t,x)\frac{\partial^2}{\partial y^2}
\eta(t,u(t,x))\frac{\partial}{\partial
x}u(t,x)^T\sigma(x)\rt]\Bigg]\d t+ \frac{\partial}{\partial
z}\alpha(u(t,x),B_t)\d B_t.
\end{aligned}
\end{equation}
We notice that
$$
\frac{\partial}{\partial x}\hat{u}(t,x)=\frac{\partial}{\partial
x}[\alpha(u(t,x),B_t)]=\lt(\frac{\partial}{\partial
y}\alpha\rt)(u(t,x),B_t)\frac{\partial}{\partial x}u(t,x),
$$
\[
\begin{aligned}
\frac{\partial^2}{\partial x^2}\hat{u}(t,x)
=\frac{\partial}{\partial x}u(t,x)\lt(\frac{\partial^2}{\partial
y^2}\alpha\rt)(u(t,x),B_t)\lt(\frac{\partial}{\partial
x}u(t,x)\rt)^T+ \lt(\frac{\partial}{\partial
y}\alpha\rt)(u(t,x),B_t)\frac{\partial^2}{\partial x^2}u(t,x).
\end{aligned}
\]
Thus, by substituting the above relations into (\ref{eq:hatu}), we
can simplify the equation (\ref{eq:hatu}) and we get
\begin{equation}
\d \hat{u}(t,x)=\lt[\mathscr{L} \hat{u}
(t,x)-f\lt(t,x,\hat{u}(t,x),\sigma(x)^T\frac{\partial}{\partial
x}\hat{u}(t,x)\rt)\rt]\d t+ g(\hat{u}(x,t))\d B_t.
\end{equation}
Consequently, $\hat{u}$ is a solution of SPDE (\ref{spde}). The
proof of the converse direction is analogous. The proof is complete.
\hfill$\cajita$

In order to summarize, we have constructed a solution of SPDE
(\ref{spde}) with the help of the fractional BDSDE
(\ref{eq:x-bdsde}), by passing through the quadratic BSDE
(\ref{eq:x-BSDE}) and the associated PDE (\ref{w-pde}) with random
coefficients:
\[
\begin{array}{cccc}
&fract.\  BDSDE \ (\ref{eq:x-bdsde})&\stackrel{Thm. \
\ref{thm:ds}}{\Longrightarrow} &quadratic \ BSDE \ (\ref{eq:x-BSDE})\\
&&&\Downarrow \\
&fract.\  SPDE \ (\ref{spde})&\stackrel{ Prop.\
\ref{prop_u}}\Longleftrightarrow &PDE \ (\ref{w-pde}) \ (with \
random\  coefficients).
\end{array}
\]

Proposition \ref{prop_u} shows that, in the case of a classical
solution, the Doss-Sussman transformation establishes a link between
PDE (\ref{w-pde}) and SPDE (\ref{spde}). This motivates us to give
the following definition of the stochastic viscosity solution. For
the cases $H\in(0,1/2]$, the reader is referred to Buckdahn and Ma
\cite{BM}, Jing and Le\'on \cite{Jl}. For the classical definition
of the viscosity solution, we refer to Crandall et al. \cite{CIL}.
\begin{definition}
A  continuous random field $\hat{u}: \Omega^{\pr}\times[0,T] \times
\R^n
 \to\R$ is called a $($stochastic$)$ viscosity solution
of equation $(\ref{spde})$ if and only if
$u(t,x)=\mathcal{E}(t,\hat{u}(t,x)), (t,x)\in[0,T]\times\R^n$ is the
viscosity  solution of equation $(\ref{w-pde})$.
\end{definition}

By following a similar argument as that developed in the proof of
Theorem 4.9 in Jing and Le\'on \cite{Jl},  our preceding discussion
leads to the following theorem:
\begin{theorem}
The stochastic field $\hat{u}: \Omega^\pr\times [0,T]\times \R^n\to
\R$ defined by $\hat{u}(t,x)\triangleq \alpha (Y_t^{t,x},B_t)$ is a
$($stochastic$)$ viscosity solution of SPDE $(\ref{spde})$.
\end{theorem}
\section*{Acknowledgements}
The author thanks Professor Rainer Buckdahn for his valuable
discussions and advice.


\begin{thebibliography}{99}
\bibitem{Amn} E.~Al\`{o}s, O.~Mazet and D.~Nualart. Stochastic calculus with respect to
Gaussian processes. {\it Annals of Probability}, \textbf{29} (2001),
766-801.
\bibitem{Be} C.~Bender. Explicit solutions of a class of linear fractional
BSDEs. \textit{Systems $\&$ Control Letters}, \textbf{54} (2005),
671-680.
\bibitem{BHOZ} F.~Biagini, Y.~Hu, B.~{\O}ksendal and
T.~Zhang. \textit{Stochastic Calculus for Fractional Brownian Motion
and Applications.} Springer, 2008.
\bibitem{bu} R.~Buckdahn. \textit{Anticipative Girsanov Transformations and
Skorohod Stochastic Differential Equations}. Memoirs of the AMS,
Vol. \textbf{111}, N.533, 1994.
\bibitem{BM} R.~Buckdahn and J.~Ma.
Stochastic viscosity solutions for nonlinear stochastic partial
differential equations. Part I. \textit{Stochastic Processes and
their Applications,} \textbf{93} (2001), 181-204.
\bibitem{CIL} M.~G.~Crandall, H.~Ishii and P-L.~Lions. User's guide to
viscosity solutions of second order partial differential equations.
\textit{Bulletin of American Mathematical Society,} \textbf{27}
(1992), 1-67.
\bibitem{Cq} L.~Coutin and Z.~Qian. Stochastic analysis, rough paths analysis
and fractional Brownian motions. {\it Probability Theory and Related
Fields}, \textbf{122} (2002), 108-140.
\bibitem{Du} L. Decreusefond and A. S. \"Ust\"unel. Stochastic analysis of
the fractional Brownian motion. \textit{Potential Analysis,}
\textbf{10} (1998), 177-214.
\bibitem{Dhp} T.~E.~Duncan, Y.~Hu and B.~Pasik-Duncan. Stochastic calculus for
fractional Brownian motion I. Theory. {\it SIAM Journal of Control
and Optimization}, \textbf{38} (2000), 582-612.
\bibitem{HP} Y.~Hu and S.~Peng. Backward stochastic differential
equation driven by fractional Brownian motion. \textit{SIAM Journal
on Control and Optimization,} {\bf 48} (2009), 1675-1700.
\bibitem{Jl} S.~Jing and J.~A.~Le\'on. Semilinear backward doubly stochastic
differential equations and SPDEs driven by fractional Brownian
motion with Hurst parameter in (0,1/2). \textit{Preprint} (2010),
available at arXiv:1005.2017.
\bibitem{Ko} M.~Kobylanski. Backward stochastic differential
equations and partial differential equations with quadratic growth.
\textit{The Annals of Probability}, Vol. \textbf{28}, No. 2, (2000),
558-602.
\bibitem{Mi} Y.~S.~Mishura. \textit{Stochastic Calculus
for Fractional Brownian Motion and Related Processes}. Springer,
2008.
\bibitem{Nr} D.~Nualart and A.~R\u{a}\c{s}canu. Differential equations
driven by fractional Brownian motion. {\it Collectanea Mathematica},
\textbf{53} (2002), 55-81.
\bibitem{PP3}  \'E.~Pardoux and S.~Peng.  Adapted solution of a backward
stochastic differential equation. {\it Systems $\&$ Control
Letters}, {\bf 14} (1990),  55--61.
\bibitem{PP2} \'E.~Pardoux and S.~Peng. Backward doubly
stochastic differential equations and systems of quasilinear
parabolic SPDEs. {\it Probability Theory and Related Fields}, {\bf
98} (1994), 209-227.
\bibitem{PP1} \'E.~Pardoux and S.~Peng. Backward
stochastic differential equations and quasilinear parabolic partial
differential equations. {\it Lecture Notes in CIS}, Vol. {\bf 176}
(1992), Springer-Verlag, 200-217.
\bibitem{Pe} S.~Peng. Probabilistic interpretation for
systems of quasilinear parabolic partial differential equations.
{\it Stochastics and Stochastics Reports}, {\bf 37} (1991), 61-74.
\bibitem{Pr} P.~Protter. {\it Stochastic Integration and Differential equations}. Springer, Berlin, 1990.
\bibitem{Rv} F.~Russo and P.~Vallois. Elements of stochastic calculus
via regularization. {\it S\'eminaire de Probabilit\'es XL, Lecture
Notes in Mathematics}, Volume {\bf1899} (2007), Springer, 147-185.
\bibitem{Rv1} F.~Russo and P.~Vallois. Forward, backward and symmetric
stochastic integration. {\it Probability Theory and Related Fields},
{\bf 97} (1993), 403-421.
\bibitem{Rv2} F.~Russo and P.~Vallois.  The generalized
covariation process and Ito-formula. {\it Stochastic Processes and
their Applications}, {\bf 59} (1995), 81-104.
\bibitem{Yo} L.~C.~Young.
An inequality of H\"older type, connected with Stieljes integration.
{\it Acta Mathematica,} \textbf{67} (1936), 251-282.
\bibitem{Za1} M.~Z\"ahle. Integration with respect to fractal
functions and stochastic calculus I. {\it Probability Theory and
Related Fields}, {\bf 111} (1998), 333-374.
\bibitem{Za2} M.~Z\"ahle. Integration with respect to fractal
functions and stochastic calculus II. {\it Mathematische
Nachrichten},  {\bf 225} (2001), 145-183.

\end{thebibliography}
\end{document}